\documentclass[a4paper,twoside,
 10pt]{amsart}
\usepackage{amssymb,amsfonts,amsmath,enumerate,color
}
\newtheorem{theor}{~~~~Theorem}

\newtheorem{cor}{~~~~Corollary}

\allowdisplaybreaks
\newtheorem{remark}{~~~~Remark}
\newtheorem{defin}{~~~~Definition}

\def\og{\leavevmode\raise.3ex\hbox{$\scriptscriptstyle\langle\!\langle$~}}
\def\fg{\leavevmode\raise.3ex\hbox{~$\!\scriptscriptstyle\,\rangle\!\rangle$}}
\begin{document}
\bibliographystyle{plain}
\title{A note on sub-Riemannian structures associated with complex Hopf fibrations}

\author
{Chengbo Li
\address{Department of
Mathematics, Tianjin University, Tianjin, 300072, China; E-mail:
chengboli@tju.edu.cn} and Huaying Zhan
\address{School of science, Tianjin university of
technology, 300384, China; Email: zhanhuaying@gamil.com}}

\begin{abstract}
Sub-Riemannian structures on odd-dimensional spheres respecting the Hopf fibration naturally appear in quantum mechanics.
We study the curvature maps for such a sub-Riemannian structure and express them using the Riemannian curvature tensor of the Fubini-Study metric of the complex projective space and the curvature form of the Hopf fibration. We also estimate the number of conjugate points of a sub-Riemannian extremal in terms of the bounds of the sectional curvature and the curvature form. It presents a typical example for the study of curvature maps and comparison theorms for a general corank 1 sub-Riemannian structure with symmetries done by C.Li and I.Zelenko in \cite{cijacobi}.
\end{abstract}

 \maketitle\markboth{Chengbo Li and Huaying Zhan}
{On sub-Riemannian structures associated with complex Hopf
fibrations}
\footnote{The authors are partially supported by TianYuan Special Funds of the National Natural Science Foundation of China (Grant No. 11126163)}
\section{Introduction}

In the present note we focus on the sub-Riemannian geodesics for sub-Riemannian structures associated with the complex dimensional Hopf fibration
$$\mathbb S^1 \hookrightarrow \mathbb S^{2n+1}\hookrightarrow \mathbb {CP}^n$$
over the complex projective space $\mathbb {CP}^n$. The motivation for the work is two-fold. On one hand, it has a natural quantum physics background. The case $n=1$ is of course the classical Hopf fibration which  was well studied and explicit formulas for geodesics
were obtained. 
However, the calculations for the high dimensional case ($n\geq 2$) become quite complicate and only partial results  were
obtained. See \cite{cmvhopf} and the references therein for  details.  On the other hand,  using the tool developed in the study of geometry of curves in Lagrange Grassmannians, we constructed in \cite{cijacobi} the curvature maps and expressed them in terms of the Riemannian  curvature
tensor of the base manifold and the curvature form of the
principle connection of the principle bundle.  However,
 the disadvantage there is the lack of examples  to be complementary for the
  theory while the sub-Riemannian structures associated with the complex Hopf fibrations can exactly play such a role. More precisely, instead of making efforts to obtain the explicit parametric expression of sub-Riemannian geodesics, we study the curvature maps of the sub-Riemannian structures associated with the complex Hopf fibration in order to have the intrinsic Jacobi equation along a sub-Riemannian extremal  so that we can establish the comparison theorems to estimate the number of conjugate points along the sub-Riemannian geodesic.

We organize the note as follows. First of all, we formulate the sub-Riemannian geodesic problems for the sub-Riemannian structures associated with complex Hopf fibrations. Secondly, we explain the constructions of the curvature maps for a contact sub-Riemannian structure and then show the expressions of the curvature maps when there are additional transverse symmetries. Finally, we apply the results to the sub-Riemannian structures associated with complex Hopf fibrations and get the comparison theorems of the estimation of the number of conjugate points along a sub-Riemannian geodesic.

\section{sub-Riemannian geodesic problem associated with complex Hopf fibrations}
We will start with a description of a sub-Riemannian structure associated with a principle $G$-connection on a principle $G$-bundle over a Riemannian manifold and then specialize to the case for the complex Hopf fibrations.

We use the standard terminology from the theory of principle $G$-bundles (see e.g. \cite{knfoundations1}). Let $\pi:P\rightarrow M$ be a principle G-bundle over a smooth manifold  $(M,g)$. For any $p\in P$ we can define in $T_pP$ the vertical subspace
$$\mathcal V_p:=\{v\in T_pP|\pi_*v=0\}$$
and $\mathcal V=\{\mathcal V_p:p\in P\}$ is usually called the \emph{vertical distribution}.  A principal $G$-connection on $P$ is a differential 1-form (connection form) on $P$ with values in the Lie algebra $\mathfrak g$  of $G$ which is G-equivariant and reproduces the Lie algebra generators of the fundamental vector fields on $P$. In other words, it is an element of $\omega\in\Omega^1(P,\mathfrak g)$  such that
\begin{itemize}
\item $\hbox{Ad}(g)(R_g^*\omega)=\omega$,where $R_g$ denotes right multiplication by $g$;
\item if $\xi\in\mathfrak g$ and $X_\xi$ is the fundamental vector field on $p$ associated to $\mathfrak g$, then $\omega(X_\xi)=\xi.$
\end{itemize}
A principle $G$-connection is equivalent to a $G$-equivariant Ehresmann connection $\mathcal H$, i.e., a smooth vector distribution $\mathcal H$ on $P$ satisfying
$$T_pP=\mathcal H_p+\mathcal V_p,\quad\quad\mathcal H_{pg}=d(R_g)_p(\mathcal H_p),\quad \forall p\in P, g\in G.$$
Such a distribution $\mathcal H$ is usually called a \emph{horizontal distribution}.

What we concerned is the case that the manifold $M$ is equipped with a Riemannian metric $\langle\cdot ,\cdot\rangle$ because a sub-Riemannian structure is then naturally associated with. Namely, the pull back $\pi^*(\langle\cdot,\cdot\rangle)$ defines an inner product on the distribution $\mathcal H$ as $\pi$ is an isomorphism between $\mathcal H_p$ and $T_{\pi(p)}M.$  The triple $(P,\mathcal H,\langle\cdot,\cdot\rangle)$ is called \emph{a sub-Riemannian structure associated with the principle $G$-bundle $\pi:P\rightarrow M$.} As a special case, the complex Hopf fibration $$\mathbb S^1 \hookrightarrow \mathbb S^{2n+1}\stackrel{\pi}{\hookrightarrow} \mathbb {CP}^n$$
is a principle $G-$bundle, where $G=U(1)\cong \mathbb S^1$ is the circle action and $\mathbb {CP}^n$ is equipped with the K\"{a}hlerian Fubini-Study metric $\langle\cdot,\cdot\rangle$.  The action of $e^{2\pi it} \in U(1)$ on $\mathbb S^{2n+1}$ is defined
by
$$e^{2\pi it}z = e^{2\pi it}(z_1, . . . , z_n) = (e^{2\pi it}z_1, . . . , e^{2\pi it}z_n).$$
We will use both real $(x_1, y_1, . . . , x_n, y_n)$ and complex coordinates $z_k = x_k + iy_k, k =
1, . . . , n$. The horizontal tangent space $H_z$ at $z\in \mathbb S^{2n+1}$ is the maximal complex subspace of the
real tangent space $T_z\mathbb S^{2n+1}$. The unit normal real vector field $N(z)$ at $z\in\mathbb S^{2n+1}$ is given by
$$N(z) =\sum_{i=1}^n x_k\partial_{x_k}+y_k\partial_{y_k}=2Re\sum_{i=1}^nz_k\partial_{z_k}.$$
The vertical real vector field
\begin{equation}\label{vertical}
V(z)=iN(z)=\sum_{i=1}^n -y_k\partial_{x_k}+x_k\partial_{y_k}=2Re \sum_{i=1}^n iz_k\partial_{z_k}
\end{equation}
is globally defined and non-vanishing  and spans the vertical distribution $\mathcal V$ of the complex
Hopf fibration $\mathbb S^1 \hookrightarrow \mathbb S^{2n+1}\stackrel{\pi}{\hookrightarrow} \mathbb {CP}^n$. A natural choice of $G$-principle connection
$\mathcal H$ is such that $\mathcal H_z$ is the orthogonal complement of $\mathcal V_z$ in $T_z\mathbb S^{2n+1}$ w.r.t. the round metric at $\forall z\in\mathbb S^{2n+1}$. The restriction of the round metric to the horizontal and vertical subspaces is denoted by $d_\mathcal H$ and $d_\mathcal V$, respectively. Concluding the above constructions we will work with a sub-Riemannian manifold that is the triple $(\mathbb S^{2n+1},\mathcal H, d_\mathcal H)$. Note that it is a special case of a sub-Riemannian structure associated with a principle bundle. Indeed, $\pi:\mathbb S^{2n+1}\rightarrow \mathbb{CP}^n$ is a Riemannian submersion, where the round metric is endowed upon $\mathbb S^{2n+1}$ and the the Fubini-Study metric is endowed upon $\mathbb{CP}^n$ (see. e.g.\cite{priemannian}), therefore $d_\mathcal H=\pi^*(\langle\cdot,\cdot\rangle)$. We finally remark that by definition the distribution $\mathcal H$ is nonholonomic or bracket-generating.

\section{Construction of curvature maps for a contact sub-Riemannian structure}
The following construction can actually be done for a general sub-Riemannian structure. However, as the sub-Riemannian structure $(\mathbb S^{2n+1},\mathcal H,d_\mathcal H)$ is of contact type, it suffices to focus on a contact sub-Riemannian structure, which also proves to be a more concrete exposition than that of the general one. For this, let  $(M, \mathcal D, \left\langle\cdot,\cdot\right\rangle)$ be a sub-Riemannian structure on $M$ and $\mathcal D$ is a contact distribution. Assume that $M$ is connected and that $\mathcal{D}$ is nonholonomic or bracket-generating. A Lipschitzian curve $\gamma: [0, T]\longrightarrow M$ is called \emph{admissible} if $\dot\gamma(t)\in \mathcal{D}_{\gamma(t)}$, for a.e. $t$. It follows from the Rashevskii-Chow theorem that any two points in $M$ can be connected by an admissible curve. One can define the length of an admissible curve $\gamma: [0, T]\longrightarrow M$ by $\int_0^T \|\dot\gamma(t)\|dt,$ where $\|\dot\gamma(t)\|=\left\langle\dot\gamma(t),\dot\gamma(t)\right\rangle^{\frac{1}{2}}.$

\subsection{Sub-Riemannian geodesics}

The length minimizing problem is to find the shortest admissible curve connecting two given points on $M$. As in Riemannian geometry, it is equivalent to the problem of minimizing the kinetic energy $\frac{1}{2}\int_0^T \|\dot\gamma(t)\|^2dt$. The problem can be regarded as an optimal control problem and its extremals can be described by the Pontryagin Maximum Principle of Optimal Control Theory (\cite{pbgmthe}). There are two different types of extremals: abnormal and normal, according to vanishing or nonvanishing of Lagrange multiplier near the functional, respectively. For the case of a contact sub-Riemannian structure, all sub-Riemannian energy (length) minimizers are the projections of normal extremals (see e.g. \cite{matour}).

Therefore we shall focus on normal extremals only. To describe them let us introduce some notations.
Let $T^*M$ be the cotangent bundle of $M$ and $\sigma$ be the canonical symplectic form on $T^*M$, i.e., $\sigma=-d\varsigma$, where $\varsigma$ is the tautological (Liouville) 1-form on $T^*M$. For each function $H:T^*M\to \mathbb R$, the Hamiltonian vector field $\vec h$ is defined by $i_{\vec h}\sigma=dh.$
Given a vector $u\in T_qM$ and a covector $p\in T_q^*M$ we denote by $p\cdot u$ the value of $p$ at $u$.
Let
\begin{equation}\label{h}
h(\lambda)\stackrel{\Delta}{=}\max_{u\in\mathcal{D}}(p\cdot u-\frac{1}{2}\|u\|^2)=\frac{1}{2}\|p|_{\mathcal{D}_q}\|^2,\quad\lambda=(p,q)\in T^*M,\ q\in M,\ p\in T^*_qM,
\end{equation}
where $p|_{\mathcal{D}_q}$ is the restriction of the linear functional $p$ to $\mathcal{D}_q$ and the norm $\|p|_{\mathcal{D}_q}\|$ is defined w.r.t. the Euclidean structure on $\mathcal D_q.$  The normal extremals are exactly the trajectories of $\dot\lambda(t)=\vec h(\lambda)$.
\subsection{Jacobi curve and conjugate points along normal extremals}
Let us fix the level set of the Hamiltonian function $h$:
$$\mathcal{H}_{c}\stackrel{\Delta}{=}\{\lambda\in T^*M| h(\lambda)=c\}, c>0$$
Let $\Pi_{\lambda}$ be the vertical subspace of
$T_{\lambda}\mathcal{H}_{c}$, i.e.
$$
\Pi_{\lambda}=\{\xi\in T_{\lambda}\mathcal{H}_c:
\pi_*(\xi)=0\},
$$
where $\pi: T^*M\longrightarrow M$ is the canonical projection.
With any normal extremal $\lambda(\cdot)$ on $\mathcal H_{c}$, one can associate a curve in a Lagrange Grassmannian which describe the dynamics of the vertical subspaces $\Pi_\lambda$ along this extremal w.r.t. the flow $e^{t\vec h}$, generated by $\vec h$. For this let
\begin{equation}\label{Jacobi}
t\longmapsto \mathfrak J_{\lambda}(t)\stackrel{\Delta} {=}e_*^{-t\vec
h}(\Pi_{e^{t\vec h}\lambda})/\{\mathbb{R}\vec h(\lambda)\}.
\end{equation}
The curve $\mathfrak J_\lambda(t)$ is the curve in the Lagrange Grassmannian of the linear symplectic space
$W_\lambda = T_\lambda\mathcal H_{c}/{\mathbb R\vec h(\lambda)}$ (endowed with the symplectic form induced in the obvious way by the canonical symplectic form $\sigma$ of $T^*M$). It is called the \emph{Jacobi curve} of the extremal $e^{t\vec h}\lambda$ (attached at the point $\lambda$).

The reason to introduce Jacobi curves is two-fold. On one hand, it can be used to construct differential invariants of sub-Riemannian structures, namely, any symplectic invariant of Jacobi curve, i.e., invariant of the action of the linear symplectic group $Sp(W_\lambda)$ on the Lagrange Grassmannian $L(W_\lambda)$, produces an invariant of the original sub-Riemannian structure. On the other hand, the Jacobi curve contains all information about conjugate points along the extremals.

Recall that time $t_0$ is called conjugate to $0$ if
\begin{equation}\label{conju}
e^{t_0\vec h}_*\Pi_\lambda\cap\Pi_{e^{t_0\vec h}\lambda}\neq 0.
\end{equation}
and the dimension of this intersection is called the multiplicity of $t_0$.
The curve  $\pi(\lambda(\cdot))|_{[0, t]}$ is $W^1_\infty$-optimal (and even $C$-optimal) if there is no conjugate point in $(0, t)$ and is not optimal otherwise. Note that (\ref{conju}) can be rewritten as: $e^{-t_0\vec h}_*\Pi_{e^{t_0\vec h}\lambda}\cap \Pi_\lambda\neq 0$, which is equivalent to $$\mathfrak J_\lambda(t_0)\cap \mathfrak J_\lambda(0)\neq 0.$$

\subsection{Curvature maps and structural equations}
For a curve $\Lambda(\cdot)$ in Lagrange Grassmannian of a linear
symplectic space $W$, satisfying very mild condition, one can
construct the complete system of symplectic
invariants(\cite{icdifferential}) and the \emph{normal moving
frame} satisfying some canonical structural equation. In
particular, for the Jacobi curve $\mathcal J_\lambda(\cdot)$, where $\lambda\in\mathcal H_{\frac{1}{2}}$,
associated with a sub-Riemannian extremal of a contact
sub-Riemannian structure, such a result reads in a simpler way.
Fix $\dim M=n.$
\begin{defin}
\label{normframe} The moving Darboux frame
$(E^\lambda_a(t),E^\lambda_b(t),E^\lambda_c(t),F^\lambda_a(t),F^\lambda_b(t),F^\lambda_{c}(t))$,
where  $$E^\lambda_a(t), E^\lambda_b(t),F^\lambda_a(t), F^\lambda_b(t)\hbox{are vectors and}\ E^\lambda_c(t)=(E^\lambda_{c_1}(t),\cdots,E^\lambda_{c_{n-3}}(t)),
F^\lambda_c(t)=(F^\lambda_{c_1}(t),\cdots,F^\lambda_{c_{n-3}}(t)),$$
 is called the normal
moving frame of $\mathcal J_\lambda(t)$, if for any $t$,
$$\mathcal J_\lambda(t)={\rm span}\{E^\lambda_a(t),E^\lambda_b(t),E^\lambda_c(t)\}$$ 
and there exists an one-parametric family of normal mappings
$(R_t(a,a),R_t(a,c),R_t(b,b),R_t(b,c),R_t(c,c))$, where $R_t(a,a),R_t(b,b)\in \mathbb R$ and $R_t(a,c),R_t(b,c)\in \mathbb R^{(n-3)\times 1}$ and $R_t(c,c)\in \mathbb R^{(n-3)\times(n-3)}$ 
is symmetric for any t, such that the moving frame
$(E^\lambda_a(t),E^\lambda_b(t),E^\lambda_c(t),F^\lambda_a(t),F^\lambda_b(t),F^\lambda_{c}(t))$,
satisfies the following structural equation:
\begin{equation}
\label{structeq}
\begin{cases}
 E_a'(t)=E_{b}(t)\\
 E_b'(t)=E_c(t)\\
E_c'(t)=F_c(t)\\
F_a'(t)=-E_a(t)R_t(a,a)-E_c(t)R_t(a,c)\\
F_b'(t)=-F_a(t)-E_b(t)R_t(b,b)-E_c(t)R_t(b,c)\\
F_c'(t)=-E_a(t)(R_t(a,c))^T-E_b(t)(R_t(b,c))^T-E_c(t)R_t(c,c).
\end{cases}
\end{equation}
\end{defin}

\begin{theor}
There exists a normal moving frame
$(E^\lambda_a(t),E^\lambda_b(t),E^\lambda_c(t),F^\lambda_a(t),F^\lambda_b(t),F^\lambda_{c}(t))$
of $\mathcal J_\lambda(t)$. Moreover, if there is another normal
moving frame $(\tilde E^\lambda_a(t),\tilde E^\lambda_b(t),\tilde
E^\lambda_c(t),\tilde F^\lambda_a(t),\tilde F^\lambda_b(t),\tilde
F^\lambda_{c}(t))$ of $\mathcal J_\lambda(t)$, then it must hold
\begin{eqnarray*}
(\tilde E^\lambda_a(t),\tilde
E^\lambda_b(t))&=&\pm
(E^\lambda_a(t),E^\lambda_b(t)),\\
\tilde E^\lambda_c(t)&=&E^\lambda_c(t)O,
\end{eqnarray*}
where $O$ is a  constant
orthornormal matrix.
\begin{remark}
If $n=3$, then $E^\lambda_c(t),F^\lambda_c(t)$ do not appear in the above construction and such a convention is understood in the remainder of the text.
\end{remark}

\end{theor}
It follows from the last theorem that there is a canonical splitting
of the subspace $\mathcal J_\lambda(t)$, i.e.,
$$\mathcal J_\lambda(t)=V_a(t)\oplus V_b(t)\oplus V_c(t),$$
where $V_a(t)=\mathbb RE_a(t),\ V_b(t)=\mathbb RE_b(t),\
V_c(t)={\rm span}\{E_c(t)\}$. Each space is endowed with the
canonical Euclidean structure, in which a or a tuple of vectors
$E_a(t),E_b(t),E_c(t)$ from some (and therefore any) normal moving
frame constitutes the orthonormal frame. For any
$s_1,s_2\in\{a,b,c\}$, the linear map from $V_{s_1}(t)$ to
$V_{s_2}(t)$ with the matrix $R_t(s_1, s_2)$ from (\ref{structeq})
in the basis $\{E_{s_1}(t)\}$ and $\{E_{s_2}(t)\}$ of $V_{s_1}(t)$
and $V_{s_2}(t)$ respectively, is independent of the choice of
normal moving frames. It will be denoted by $\mathfrak R_t(s_1,
s_2)$ and it is called the \emph {$(s_1,s_2)$-curvature map of the
curve $\Lambda(\cdot)$ at time $t$}.
\subsection{Expressions of the  curvature maps and the comparison theorems}
The construction above helps to find very fruitful additional structures in the cotangent bundle $T^*M$.
The structural equation (\ref{structeq}) for the Jacobi curve $\mathcal J_\lambda(t)$ can be seen as the intrinsic Jacobi equation along the extremal $e^{t\vec h}\lambda$ and the curvature maps are the coefficients of this Jacobi equation.

Since there is a canonical splitting of $\mathcal J_\lambda(t)$ and taking into account that $\mathcal J_\lambda(0)$ and $\Pi_\lambda$ can be naturally identified, we have the canonical splitting of $\Pi_\lambda$:
$$\Pi_\lambda=\mathcal V_a(\lambda)\oplus \mathcal V_b(\lambda)\oplus \mathcal V_c(\lambda),\dim \mathcal V_a(\lambda)=\dim \mathcal V_b(\lambda)=1,\ \dim\mathcal V_c(\lambda)=n-3,$$
where $\mathcal V_s(\lambda)=V_s(0), s=a,b,c$. Moreover, let $\mathfrak R_\lambda(s_1,s_2): \mathcal V_{s_1}(\lambda)\rightarrow \mathcal V_{s_2}(\lambda)$  and the $\mathfrak R_\lambda:\Pi_\lambda\rightarrow\Pi_\lambda$ be the $(s_1,s_2)$-curvature map.
These maps are intrinsically related to the sub-Riemannian structure.
They are called the \emph{$(s_1,s_2)$-curvature}.

In the  Riemannian case, the  curvature map is expressed in terms of Riemannian curvature tensor and the structural equations are actually the Jacobi equations in Riemannian geometry. For a sub-Riemannian structures $(P,\mathcal H,\langle\cdot,\cdot\rangle)$ associated with a principle connection $\mathcal H$ on a $G$-bundle $\pi:P\rightarrow M$ with one dimensional fibers, it turns out that the big curvature map is a combination of Riemannian curvature tensor of $M$ and the curvature form. To be more precise, let $\omega$ be the connection 1-form of $\mathcal H$ then $d\omega$ is the curvature form and  it induces a 1-1 tensor on $M$
$$g(JX,Y)=d\omega(X,Y).$$
 A general formula of the curvature maps using the Riemannian curvature tensor on $M$ and the tensor $J$, together with their covariant derivatives can be found in \cite{cijacobi}. 

Now let us specialize to the case of a sub-Riemannian structure $(\mathbb S^{2n+1},\mathcal H, d_\mathcal H)$. It can be shown that $(J,g)$ defines a K\"{a}hlerian structures on $\mathbb{CP}^n$. See e.g. \cite[Chapter 3]{priemannian}. In this case the curvature maps read in a very simple form.
For this, let us first of all give more explicit description of the subspaces $\mathcal
V_z(\lambda)$. As the tangent space of the fibers of $T^*\mathbb S^{2n+1}$ can be naturally identified with the fibers themselves (the fibers are linear spaces), one can show that,
 $$\mathcal V_a(\lambda)=\mathcal{\mathcal H}_{z}^\bot,$$
 where $\mathcal{\mathcal H}_{z}^\bot$ is the annihilator of $\mathcal H$, namely,
 $$\mathcal{H}_{z}^\bot=\{p\in T^*_z\mathbb S^{2n+1}:p\cdot v=0,\ \forall v\in \mathcal H_z\}.$$
 Since the
 Moreover, we have that
\begin{equation}\label{ident3}
\mathcal V_b(\lambda)\oplus \mathcal V_c(\lambda)
\sim\mathcal {H}_z^*\sim \mathcal{H}_z.
\end{equation}
Since the $\pi_*:\mathcal H_{z}\rightarrow T_{\pi(z)} \mathbb C\mathbb P^{n}$ is an isometry for all $z\in \mathbb S^{2n+1}$, we also take $J_z$ as a antisymmetric operator on $\mathcal H_z$. So, under the above identifications, one can show that

\begin{equation}\label{ident4}
\mathcal V_b(\lambda)=\mathbb{R}J_zp,\quad \mathcal V_c(\lambda)=(\mathbb
R J_z p)^\bot.
\end{equation}

Actually, $d\omega$ can be seen as a magnetic field and $J$ can be seen as a Lorenzian force on $\mathbb C\mathbb P^n$. The projection by $\pi$ of all sub-Riemannian geodesics describes all possible motion of a charged particle (with any possible charge) given by the magnetic field $d\omega$ on the Riemannian manifold $\mathbb {CP}^n$(see e.g. \cite[Chapter 12]{matour} and the references therein). 


Define $u_0: T^*\mathbb S^{2n+1}\to \mathbb R$ by
$u_0(p,z):=p\cdot V(z),\ z\in \mathbb S^{2n+1},p\in T^*_z\mathbb S^{2n+1},$
where $V$ is the vertical vector field defined in \eqref{vertical}.
 As before $\lambda=(p,z)\in\mathcal H_{\frac{1}{2}},z\in \mathbb S^{2n+1},p\in T_z^*\mathbb S^{2n+1}$,  any $v\in T_\lambda T^*_z \mathbb S^{2n+1} (\sim T^*_z\mathbb S^{2n+1}\sim T_z \mathbb S^{2n+1})$, we have a vector $v^h:=\pi_*v\in T_{\pi(z)}\mathbb {CP}^n$; conversely, given any $X\in T_{\pi (z)}\mathbb{CP}^n$, there is a unique $X^v\in \mathcal H_z\subset T_z\mathbb{S}^{2n+1} (\sim T_{(p,z)} T_z \mathbb{S}^{2n+1}), \ p\in T_z^*\mathbb{S}^{2n+1}$.

\begin{theor}
\label{Kahler}
Let $(\langle\cdot,\cdot\rangle,J)$ be the K\"{a}hlerian structure on $\mathbb {CP}^n$ and $R^\nabla$ the Riemannian curvature tensor of $\langle\cdot,\cdot\rangle$. Then for $\forall v\in \mathcal V_c(\lambda)$,
\begin{eqnarray*}
 g((\mathfrak{R}_\lambda(c,c)(v))^h, v^h)&=&g(R^\nabla(p^h, v^h)p^h,v^h)+\frac{u_0^2}{4}\|v\|^2,\\
\mathfrak R_\lambda(b,c)(v)&=&g(R^\nabla(p^h, Jp^h)p^h,v^h)\mathcal E_b(\lambda),\\
\rho_\lambda(b,b)&=&g(R^\nabla(p^h, Jp^h)p^h,Jp^h)+u_0^2,\\
\mathfrak R_\lambda(c,a)&=&0\quad \hbox{and}\quad \mathfrak R_\lambda(a,a)=0,
\end{eqnarray*}
where $\mathcal E_b(\lambda)$ and $\rho_\lambda(b,b)$ are defined by
$$\mathcal E_b(\lambda)=(Jp^h)^v,\quad \mathfrak R_\lambda(b,b)v_b=\rho_\lambda(b,b)v_b,\ \forall v_b\in \mathcal V_b(\lambda).$$
\end{theor}

 Now let us recall an estimate on the bounds of sectional curvature of Fubini-Study metric on $\mathbb{CP}^n$ from the theory of Riemannian submersion (see e.g. \cite[chapter 3]{priemannian})

\begin{theor}\label{Oneill}
Let $\sec(g)$ be the Riemannian sectional curvature of $g$ on $\mathbb{CP}^n$. Then  $\sec(g)\in[1,4].$ Moreover, the estimate of the bounds is sharp, namely, the values 1 and 4 are achieved.
\end{theor}

Using the tool of the Generalized Sturm Theorem for curves in Lagrangian Grassmannians,  we obtained the comparison theorems of estimation of number of conjugate points along a  sub-Riemannian extremal of a contact sub-Riemannian structure with symmetries and satisfying some compatible condition (see \cite{cijacobi} and the references therein). 
\begin{theor}
The number of conjugate points $\sharp_T\bigl(\lambda(\cdot)\bigr)$ to 0 on $(0,T]$ along $\lambda(\cdot)$ satisfies the following inequality
\begin{equation}
\label{conjest}
Z_T(1+\bar u_0^2, 1+\frac{1}{4}\bar u_0^2)\leq \sharp_T(\lambda(\cdot))\leq   Z_T(4+\bar u_0^2, 4+\frac{1}{4}\bar u_0^2),
\end{equation}
where
$$Z_T(\omega_b,\omega_c)=(n-3)[\frac{T\sqrt{\omega_c}}{\pi}]+[\frac{T\sqrt{\omega_b}}{2\pi}]+
\sharp_T\{\tan(\frac{\sqrt{\omega_b}}{2}x)-\frac{\sqrt{\omega_b}}{2}x=0\}.$$
\end{theor}
\begin{cor}\label{estimation}
Under the same estimates  the following statement hold for a normal sub-Riemannian extremal on $\mathcal H_{\frac{1}{2}}\cap\{u_0=\bar u_0\}$:

\begin{enumerate}
 \item There is no conjugate points to $0$ in the interval $\bigl(0, \frac{\pi}{\sqrt{4+\frac{1}{4}\bar u_0^2}}\bigr)$;
 \item There is at least $(n-3)$ conjugate points to $0$ in the interval $\bigl(0,\frac{2\pi}{\sqrt{4+\bar u_0^2}}
 \bigr]$ and there are at least $(n-2)$ in the interval  $\bigl(0,\frac{2\pi}{\sqrt{1+\bar u_0^2}}
 \bigr]$
\end{enumerate}
\end{cor}

\setlength\parindent{0pt}$\bullet$ {\bf  Relation to quantum systems}
The sub-Riemannian minimization problem for $(\mathbb S^{2n+1},\mathcal H, \langle\cdot,\cdot\rangle)$  , the initial and end points represent
initial and target states of the system and to find a minimizer is equivalent to find a path which transfer the minimal energy from the initial to the target states (see \cite{cmvhopf}). So, we have shown that a sub-Riemannian geodesic $\gamma(\cdot)=\pi(\lambda(\cdot))$ always transfers the minimum energy from the state $\gamma(0)$ to the state $\gamma(T),\ \forall T<\frac{\pi}{\sqrt{4+\frac{1}{4}\bar u_0^2}}$ but fails to do this from the state $\gamma(0)$ to the state $\gamma(T),\ \forall T\geq \frac{2\pi}{\sqrt{4+\bar u_0^2}}$.
\bibliography{myrefe}

\begin{thebibliography}{1}

\bibitem{cmvhopf}
D.~C. Chang, I.~Markina, and A.~Vasil¡¯ev.
\newblock Hopf fibration: Geodesics and distances.
\newblock {\em J. Geom. Phys.}, 61(6):986--1000, 2011.

\bibitem{knfoundations1}
S.~Kobayashi and K.~Nomizu.
\newblock {\em Foundations of Differential Geometry}, volume~1.
\newblock Wiley-Interscience, 1 edition, 1996.

\bibitem{cijacobi}
C.~Li and I.~Zelenko.
\newblock Jacobi equations and comparison theorems for corank 1
  sub-{R}iemannian structures with symmetries.
\newblock {\em J. Geom. Phys.}, 61:781--807, 2011.

\bibitem{matour}
R.~Montgomery.
\newblock {\em A Tour of Subriemannian Geometries, Their Geodesics, and
  Applications}.
\newblock Mathematical Surveys and Monographs, Volume 91. American Mathematical
  Society, 2002.

\bibitem{priemannian}
P.~Petersen.
\newblock {\em Riemannian Geometry}.
\newblock Springer-Verlag, 2 edition, 2006.

\bibitem{pbgmthe}
L.~S. Pontryagin, V.~G. Boltyanskii, R.~V. Gamkrelidze, and E.~F. Mischenko.
\newblock {\em The Mathematical Theory of Optimal Processes}.
\newblock Wiley, New York, 1962.

\bibitem{icdifferential}
I.~Zelenko and C.~Li.
\newblock Differential geometry of curves in {L}agrange {G}rassmannians with
  given {Y}oung diagram.
\newblock {\em Differ. Geom. Appl.}, 27(6):723--742, 2009.

\end{thebibliography}
\end{document}